\documentclass[a4paper,english,reqno]{amsart}
\usepackage[T1]{fontenc}
\usepackage[latin1]{inputenc}
\usepackage{xcolor}
\usepackage{babel}
\usepackage{amsthm}
\usepackage{amstext}
\usepackage{graphicx}
\usepackage{amssymb}
\usepackage{esint}
\PassOptionsToPackage{normalem}{ulem}
 \makeatletter
\numberwithin{equation}{section}
\numberwithin{figure}{section}
\theoremstyle{plain}
\newtheorem{thm}{Theorem}[section]
  \theoremstyle{definition}
  \newtheorem{defn}[thm]{Definition}
 \theoremstyle{definition}
  \newtheorem{example}[thm]{Example}
  \theoremstyle{plain}
  \newtheorem{cor}[thm]{Corollary}
  \theoremstyle{plain}
  \newtheorem{prop}[thm]{Proposition}
  \theoremstyle{plain}
  \newtheorem{lem}[thm]{Lemma}
  \theoremstyle{definition}
  \newtheorem{condition}[thm]{Condition}
  \theoremstyle{remark}
  \newtheorem{rem}[thm]{Remark}
  \theoremstyle{plain}
  \newtheorem{assumption}[thm]{Assumption}
  \theoremstyle{remark}
  \newtheorem*{acknowledgement*}{Acknowledgement}

\usepackage{mathptmx}
\usepackage{bbm}
\usepackage{bm}
 \newcommand{\e}{\mathrm{e}}
\newcommand{\1}{\mathbbm{1}}
\newcommand{\N}{\mathbb{N}}

\newcommand{\R}{\mathbb{R}}

\renewcommand{\P}{\mathcal{P}}
\renewcommand{\Pi}{\pi}
\renewcommand{\phi}{\varphi}
\renewcommand{\tilde}{\widetilde}
\renewcommand{\emptyset}{\varnothing}

\renewcommand{\lg}{\log}

\newcommand\diam{\mathrm{diam}}

\DeclareMathOperator*{\Int}{Int}

\DeclareMathOperator*{\card}{card}

\DeclareMathOperator*{\Con}{Con}

\DeclareMathOperator*{\dom}{dom}

 \makeatother

\begin{document}
\selectlanguage{english}

\title[Regularity of the multifractal spectrum of conformal IFS]{Regularity of multifractal spectra of conformal iterated function
systems}

\author{Johannes Jaerisch and Marc Kesseb\"ohmer}

\date{\today}

\subjclass[2000]{Primary 37C45; Secondary 37D45, 37D35}
\begin{abstract}
We investigate multifractal regularity for infinite conformal iterated
function systems (cIFS). That is we determine to what extent the multifractal
spectrum depends continuously on the cIFS and its thermodynamic potential.
For this we introduce the notion of regular convergence for families
of cIFS not necessarily sharing the same index set, which guarantees
the convergence of the multifractal spectra on the interior of their
domain. In particular, we obtain an Exhausting Principle for infinite
cIFS allowing us to carry over results for finite to infinite systems,
and in this way to establish a multifractal analysis without the usual
regularity conditions. Finally, we discuss the connections to the
$\lambda$-topology introduced by Roy and Urba\'{n}ski.
\end{abstract}

\address{{AG Dynamical Systems and Geometry},
FB-3 Mathematik und Informatik, Universit\"at Bremen, Bibliothekstrasse
1, 28359 Bremen, Germany.}

\email{mhk@math.uni-bremen.de}

\email{jogy@math.uni-bremen.de}

\maketitle

\section{Introduction and statement of results}

The theory of multifractals has its origin at the boarderline between
statistical physics and mathematics - classical references are e.~g.
\cite{FrischParisi:85,Mandelbrot:74,Mandelbrot:88,HalseyJensenKadanoff:86}.
In this paper we study multifractal spectra in the setting of infinite
conformal iterated functions systems (cIFS). These systems are given
by at most countable families $\Phi=\left(\phi_{e}:X\to X\right)_{e\in I}$
, $I\subset\N$, of conformal contractions on a compact connected
subset $X$ of the euclidean space $\left(\R^{D},\Vert\cdot\Vert\right)$,
$D\geq1$. The set of cIFS with fixed phase space $X$ will be denoted
by $\mathrm{CIFS}\left(X\right)$ (see Section \ref{sec:Preliminaries-from-CIFS}
for definitions). For $\omega\in I^{\N}$ we let $\omega_{|k}:=\omega_{1}\cdots\omega_{k}$
and $\phi_{\omega_{|k}}:=\phi_{\omega_{1}}\circ\cdots\circ\phi_{\omega_{k}}$.
Then for each $\omega\in I^{\N}$ the intersection $\bigcap_{k=1}^{\infty}\phi_{\omega_{|k}}\left(X\right)$
is always a singleton given rise to a canonical coding map $\pi_{\Phi}:I^{\N}\to X$.
Its image $\Lambda_{\Phi}:=\pi_{\Phi}\left(I^{\N}\right)$ will be
called the \emph{limit set} of $\Phi$. Given a H\"older continuous
function $\psi:I^{\N}\longrightarrow\R$ the multifractal analysis
of the system $\Phi$ with respect to the potential $\psi$ is in
our context understood to be the analysis of the level sets\[
\mathcal{F}_{\alpha}:=\pi_{\Phi}\left\{ \omega\in I^{\N}:\lim_{k\rightarrow\infty}\frac{S_{k}\psi\left(\omega\right)}{\log\bigl\Vert\phi_{\omega_{|k}}'\bigr\Vert_{X}}=\alpha\right\} \]
in terms of their Hausdorff dimension $f\left(\alpha\right):=\dim_{H}\left(\mathcal{F}_{\alpha}\right)$.
In here, $S_{k}\psi:=\sum_{n=0}^{k-1}\psi\circ\sigma^{n}$ denotes
the Birkhoff sum of $\psi$ with respect to the shift map $\sigma:I^{\N}\to I^{\N}$
on the symbolic space, and $\bigl\Vert\phi_{\omega_{|k}}'\bigr\Vert_{X}:=\sup_{x\in X}\bigl|\phi_{\omega_{|k}}'\left(x\right)\bigr|$
with $\bigl|\phi_{\omega_{|k}}'\left(x\right)\bigr|$ denoting the
operator norm of the derivative. A good reference for this kind of
multifractal analysis is provided e.~g. in \cite{pesindimensiontheoryMR1489237}.

Let us define the \emph{geometric potential function} associated with
$\Phi$ by $\zeta:I^{\N}\rightarrow\R_{0}^{-}$, $\zeta\left(\omega\right):=\log\left|\phi'_{\omega_{1}}\left(\pi\left(\sigma\left(\omega\right)\right)\right)\right|$.
It is well known that in the case of finite cIFS, that is $\card\left(I\right)<\infty$,
$f$ can be related to the Legendre transform of the \emph{free energy
function} $t:\R\to\R$, which is defined implicitly by the pressure
equation (cf. Definition \ref{def:The-topological-pressure}) \begin{equation}
\P\left(t\left(\beta\right)\zeta+\beta\psi\right)=0,\quad\beta\in\R.\label{eq:pressure equation}\end{equation}
More precisely, there exists a closed finite interval $J\subset\R$
such that for all $\alpha\in J$ we have \begin{equation}
f\left(\alpha\right)=-t^{*}\left(-\alpha\right):=-\sup_{\beta}\left\{ -\beta\alpha-t\left(\beta\right)\right\} =\inf_{\beta}\left\{ t\left(\beta\right)+\beta\alpha\right\} ,\label{eq:multifractal sp as legendre transform}\end{equation}
and for $\alpha\notin J$ we have $\mathcal{F}_{\alpha}=\emptyset$
(\cite[Theorem 21.1]{pesindimensiontheoryMR1489237}, \cite{Schmeling:99}).
If we consider infinite cIFS, i.~e. $\card\left(I\right)=\card\left(\N\right)$,
we have to take into account that the pressure function might behave
irregularly and hence it is not always possible to find a solution
of (\ref{eq:pressure equation}). For the special case in which (\ref{eq:pressure equation})
has a unique solution the multifractal analysis has been discussed
in \cite[Section 4.9]{urbanskimauldin-gdmsMR2003772}. Further interesting
results on the spectrum of local dimension for Gibbs states can be
found \cite{RoyUrbanski:08}. 

Our first task is to generalise this concept to the case when the
free energy cannot be defined by the unique solution of (\ref{eq:pressure equation}).
This leads to the following modified definition of the free energy
function. 
\begin{defn}
\label{def:temperature functions}Let $\Phi\in\mathrm{CIFS}\left(X\right)$
and $\psi:I^{\N}\longrightarrow\R$ be a potential function. Then
the \emph{free energy function} $t:\R\longrightarrow\R\cup\left\{ \infty\right\} $
for the pair $\left(\Phi,\psi\right)$ is given by\begin{equation}
t\left(\beta\right):=\inf\left\{ t\in\R:\P\left(t\zeta+\beta\psi\right)\leq0\right\} .\label{eq:FreeEnergy}\end{equation}
 
\end{defn}
Notice that our definition of the free energy function generalises
the definition given for the multifractal analysis presented in \cite[Section 4.9]{urbanskimauldin-gdmsMR2003772}
or in \cite{MR2343687}, where the existence of a zero of the pressure
function $t\longmapsto\P(t\zeta+\beta\psi)$ is always required. Our
definition is rather in the spirit of \cite[Theorem 4.2.13]{urbanskimauldin-gdmsMR2003772},
which gives a version of Bowen's formula, without assuming a zero
of the pressure function to exist. More precisely, we have \[
\dim_{H}\left(\Lambda_{\Phi}\right)=\inf\left\{ t\in\R:\P\left(t\zeta\right)\leq0\right\} ,\]
which immediately implies that $t\left(0\right)=\dim_{H}\left(\Lambda_{\Phi}\right)$.
In fact, Lemma \ref{lem:ProperConvex} shows that Definition \ref{def:temperature functions}
gives rise to a proper convex function. This concept of the free energy
function has been investigated further in \cite{JaerischKessLamei:10}
as a special case of the \emph{induced topological pressure} for arbitrary
countable Markov shifts. We would like to point out that this new
formalism gives rise to further interesting exhausting principles
similar to Example \ref{exa:ExhaustionPrinciple} and Corollary \ref{cor:-ExhaustingPrinciple}
above. 

To state our first main result we set \[
\alpha_{-}:=\inf\left\{ -t^{-}\left(x\right):x\in\Int\left(\dom\left(t\right)\right)\right\} \;\mbox{and }\;\alpha_{+}:=\sup\left\{ -t^{+}\left(x\right):x\in\Int\left(\dom\left(t\right)\right)\right\} ,\]
where $t^{+}$, resp. $t^{-}$, denotes the derivative of $t$ from
the right, resp. from the left, $\Int\left(A\right)$ denotes the
interior of the set $A$, and $\dom\left(t\right):=\left\{ x\in\R:t\left(x\right)<+\infty\right\} $
refers to the effective domain of $t$. 
\begin{thm}
\label{thm:tdach als hausdorff dim} For $\alpha\in\R$ we have $f\left(\alpha\right)\leq\max\left\{ -t^{*}\left(-\alpha\right),0\right\} $
and for $\alpha\in\left(\alpha_{-},\alpha_{+}\right)$ we have $f\left(\alpha\right)=-t^{*}\left(-\alpha\right)$. 
\end{thm}
This first main result is essentially a consequence of the multifractal
regularity property of sequences of tuples $\left(\Phi^{n},\psi^{n}\right)_{n}$
of iterated function systems and potentials, which is the second main
concern of this paper.

We adapt the definition of pointwise convergence in $\mathrm{CIFS}\left(X\right)$
as used by Roy and Urba\'{n}ski in \cite{MR2183304} to our setting,
allowing us to investigate also families of cIFS with associated potentials
not sharing the same index set $\N$. To simplify notation let us
write $\Vert h\Vert_{\Omega}:=\sup_{\omega\in\Omega}\left\vert h\left(\omega\right)\right\vert $
for the supremum norm of the map $h:\Omega\longrightarrow\left(V,\vert\cdot\vert\right)$
from $\Omega$ to the normed space $\left(V,\vert\cdot\vert\right)$.
For $\Phi^{1},\Phi^{2}\in\mathrm{CIFS}\left(X\right)$ we define \begin{equation}
\rho\left(\Phi^{1},\Phi^{2}\right):=\sum_{i\in I_{1}\cap I_{2}}2^{-i}\left(\bigl\Vert\phi_{i}^{1}-\phi_{i}^{2}\bigr\Vert_{X}+\bigl\Vert(\phi_{i}^{1})'-(\phi_{i}^{2})'\bigr\Vert_{X}\right)+\sum_{i\in I_{1}\triangle I_{2}}2^{-i},\label{eq:RhoMetrik}\end{equation}
where $A\triangle B$ denotes the usual symmetric difference of the
sets $A$ and $B$. It will turn out that $\rho$ defines a metric
on $\mathrm{CIFS}\left(X\right)$. For $\omega\in\N^{k}$ and $k\in\N$
we let $\left[\omega\right]:=\left\{ \tau\in\N^{\N}:\tau_{|k}=\omega\right\} $
denote the \emph{cylinder set} of $\omega$. 

In order to set up a multifractal spectrum we restrict our analysis
to families of H\"older continuous functions $\psi^{n}:I_{n}^{\N}\longrightarrow\R$
and $\psi:I^{\N}\longrightarrow\R$ with $I_{n}\subset I\subset\N$,
$n\in\N$. 
\begin{defn}
\label{def:pointwise convergence}\renewcommand{\theenumi}{\Alph{enumi}}
We say that $\left(\Phi^{n},\psi^{n}\right)_{n}\longrightarrow\left(\Phi,\psi\right)$
\emph{converges pointwise if}, 
\begin{enumerate}
\item $\Phi^{n}\longrightarrow\Phi$ in the $\rho$-metric and
\item \label{enu:locallyUniformPsi} for all $k\in I$ we have~ ${\displaystyle \lim_{n\to\infty}\bigl\Vert\psi^{n}-\psi\bigr\Vert_{\left[k\right]\cap I_{n}^{\N}}=0}$.
\end{enumerate}
\end{defn}
Notice, that the convergence in $\rho$-metric implies that $\bigl\Vert\psi^{n}-\psi\bigr\Vert_{\left[k\right]\cap I_{n}^{\N}}$
in (\ref{enu:locallyUniformPsi}) is well defined for all sufficiently
large $n$. For a further discussion of the above defined property
see also the remark succeeding Lemma \ref{lem:ConditionM*forZeta+Psi}. 

As discussed in \cite{MR2183304} pointwise convergence topology leads
to discontinuities of the Hausdorff dimension of the limit sets. By
introducing a weaker topology called the $\lambda$-topology in \cite{MR2183304}
the Hausdorff dimension of the limit set depends continuously on the
system (see also \cite{lambdatopologypreprint}). Convergence in $\lambda$-topology
requires the additional condition (\ref{eq:lambda topology property})
below. As a corollary we will also establish the continuity of the
Hausdorff dimension under weaker assumptions. 

We are going to employ similar assumptions on the convergence of the
pairs $\left(\Phi^{n},\psi^{n}\right)_{n}$ and $\left(\Phi,\psi\right)$
to obtain continuity of the multifractal spectra. This is the purpose
of the following definition. For this let $\zeta^{n}$ denote the
geometric potential associated with $\Phi^{n}$.
\begin{defn}
\label{def:regularconvergence-1} We say that$\left(\Phi^{n},\psi^{n}\right)_{n}$
\emph{converges regularly} to $\left(\Phi,\psi\right)$, if $\left(\Phi^{n},\psi^{n}\right)_{n}\longrightarrow\left(\Phi,\psi\right)$
converges pointwise, and if for $t,\beta\in\R$ with $\P(t\zeta+\beta\psi)<\infty$
there exists $k\in\N$ and a constant $C>0$ such that for all $n\in\N$
and all $\omega\in\left(I_{n}\right)^{k}$ we have\[
\exp\sup_{\tau\in I_{n}^{\N}\cap\left[\omega\right]}\left(S_{k}\left(t\zeta^{n}+\beta\psi^{n}\right)\left(\tau\right)\right)\le C\exp\sup_{\rho\in I^{\N}\cap\left[\omega\right]}\left(S_{k}\left(t\zeta+\beta\psi\right)\left(\rho\right)\right).\]

\end{defn}
The assumption in Definition \ref{def:regularconvergence-1} is similar
to the corresponding inequality in the definition of the $\lambda$-topology
in \cite{MR2183304} but depends additionally on the potentials $\psi^{n}$
and $\psi$. For particular cases we will show that the convergence
$\Phi^{n}\longrightarrow\Phi$ in the $\lambda$-topology immediately
implies the conditions in Definition \ref{def:regularconvergence-1}.
This is demonstrated in the following example providing an analysis
of the (inverse) Lyapunov spectrum. This example is covered by Proposition
\ref{pro:connection with lambda topology} (\ref{enu:condition1})
stated in Section \ref{sec:Connection with lambda}. 
\begin{example}
[${\bm \lambda}$\bf{-topology}] \label{exa:lambda topology} Let
$\Phi^{n}=\left(\phi_{e}^{n}\right)_{e\in I_{n}}$, $\Phi=\left(\phi_{e}\right)_{e\in\N}$
be elements of $\mathrm{CIFS}\left(X\right)$ with $\Phi^{n}\rightarrow\Phi$
converging in the $\lambda$-topology and let $\psi^{n}=\psi=1$.
Then $\left(\Phi^{n},\psi^{n}\right)_{n}\longrightarrow\left(\Phi,\psi\right)$
converges regularly.
\end{example}
The second example -- eventhough straightforward to verify -- is not
only interesting for itself but will be of systematic importance for
the proof of Theorem \ref{thm:tdach als hausdorff dim}. See also
Remark \ref{rem:ExhaustionPrinciple} and Example \ref{cor:-ExhaustingPrinciple}
for further discussion of this example.
\begin{example}
[\bf {Exhausting Principle I}] \label{exa:ExhaustionPrinciple}

Let $\Phi=\left(\phi_{e}\right)_{e\in\N}$ be an element of $\mathrm{CIFS}\left(X\right)$
and $\psi:I^{\N}\longrightarrow\R$ be H\"older continuous. Define $I_{n}:=I\cap\left\{ 1,\dots n\right\} $,
$n\in\N$ and let $\Phi^{n}=\left(\phi_{e}\right)_{e\in I_{n}}$ and
$\psi^{n}:=\psi\big|_{I_{n}^{\N}}$. Then $\left(\Phi^{n},\psi^{n}\right)_{n}\longrightarrow\left(\Phi,\psi\right)$
converges regularly.
\end{example}
If the multifractal regularity property is satisfied we are able to
prove the regularity of the free energy functions.
\begin{thm}
\label{thm:convergence of temperature function}If $\left(\Phi^{n},\psi^{n}\right)_{n}\longrightarrow\left(\Phi,\psi\right)$
converges regularly then $t_{n}$ converges pointwise to $t$ on $\R$. 
\end{thm}
To state our second main result on the regularity of the multifractal
spectra let \[
\mathcal{F}_{\alpha}^{n}:=\pi_{\Phi^{n}}\left\{ \omega\in I_{n}^{\N}:\lim_{k\rightarrow\infty}\frac{S_{k}\left(\psi^{n}\right)\left(\omega\right)}{\log\bigl\Vert(\phi_{\omega_{|k}}^{n})'\bigr\Vert_{X}}=\alpha\right\} ,\quad f_{n}\left(\alpha\right):=\dim_{H}\left(\mathcal{F}_{\alpha}^{n}\right),\]
and with $t_{n}$ denoting the free energy function of $\left(\Phi^{n},\psi^{n}\right)$
let \[
\alpha_{-}^{n}:=\inf\left\{ -t_{n}^{-}\left(x\right):x\in\Int\left(\dom\left(t_{n}\right)\right)\right\} \;\mbox{and }\;\alpha_{+}^{n}:=\sup\left\{ -t_{n}^{+}\left(x\right):x\in\Int\left(\dom\left(t_{n}\right)\right)\right\} .\]

\begin{thm}
\label{thm:multifractal regularity}Let $\Phi^{n}=\left(\phi_{e}^{n}\right)_{e\in I_{n}}$,
$\Phi=\left(\phi_{e}\right)_{e\in I}$ be elements of $\mathrm{CIFS}\left(X\right)$
and $\psi^{n},\psi$ be H\"older potentials such that $\left(\Phi^{n},\psi^{n}\right)_{n}\longrightarrow\left(\Phi,\psi\right)$
converges regularly. Then for each $\alpha\in\left(\alpha_{-},\alpha_{+}\right)$
we have 
\begin{itemize}
\item \textup{$\lim_{n\to\infty}-t_{n}^{*}\left(-\alpha\right)=f\left(\alpha\right)=-t^{*}\left(-\alpha\right)$,}
\item $f_{n}\left(\alpha\right)=-t_{n}^{*}\left(-\alpha\right)$, for all
$n$ sufficiently large. 
\end{itemize}
In particular, we have $\limsup_{n}\alpha_{-}^{n}\le\alpha_{-}\le\alpha_{+}\le\liminf_{n}\alpha_{+}^{n}.$
If additionally $\sup\dom\left(t\right)=+\infty$ then $\liminf_{n}\alpha_{-}^{n}\ge\alpha_{-}$,
whereas, if $\inf\dom\left(t\right)=-\infty$ then $\limsup_{n}\alpha_{+}^{n}\le\alpha_{+}$. 
\end{thm}
Combining the above theorems with Example \ref{exa:ExhaustionPrinciple}
we obtain the following application of our analysis. 
\begin{cor}
[\bf {Exhausting Principle II}] \label{cor:-ExhaustingPrinciple}
Let $\Phi=\left(\phi_{e}\right)_{e\in\N}$ be an element of $\mathrm{CIFS}\left(X\right)$,
$\psi:I^{\N}\longrightarrow\R$ be H\"older continuous, and \textup{$\Phi^{n}=\left(\phi_{e}\right)_{e\in I_{n}}$
and $\psi^{n}:=\psi\big|_{I_{n}^{\N}}$ with }$I_{n}:=I\cap\left\{ 1,\dots n\right\} $,
$n\in\N$. Then for each $\alpha\in\left(\alpha_{-},\alpha_{+}\right)$
we have \textup{$\lim_{n}-t_{n}^{*}\left(-\alpha\right)=f\left(\alpha\right)=-t^{*}\left(-\alpha\right)$
and }$f_{n}\left(\alpha\right)=-t_{n}^{*}\left(-\alpha\right)$, for
all $n$ sufficiently large. For the boundary points of the spectrum
we have the following.
\begin{enumerate}
\item If $\dom\left(t\right)=\R$ then $\lim_{n\to\infty}\alpha_{\pm}^{n}=\alpha_{\pm}$,
\item if $\sup\dom\left(t\right)<+\infty$ then $\limsup_{n}\alpha_{-}^{n}=-\infty$
and for all $\alpha<\alpha_{-}$ we have \[
\limsup_{n\to\infty}-t_{n}^{*}\left(-\alpha\right)\leq f\left(\alpha\right),\]

\item \label{enu:LowerBoundf}if $\inf\dom\left(t\right)>-\infty$ then
$\liminf_{n}\alpha_{+}^{n}=+\infty$ and for all $\alpha>\alpha_{+}$
we have \[
\limsup_{n\to\infty}-t_{n}^{*}\left(-\alpha\right)\leq f\left(\alpha\right).\]
 
\end{enumerate}
\end{cor}
In Example \ref{exa:-LuerothGeneralized} below we demonstrate how
the lower bound on $f$ stated in Corollary \ref{cor:-ExhaustingPrinciple}
(\ref{enu:LowerBoundf}) can be applied. 

Note that by virtue of Proposition \ref{pro:connection with lambda topology}
we have on the one hand that the convergence $\Phi^{n}\longrightarrow\Phi$
in the $\lambda$-topology implies that $\left(\Phi^{n},0\right)_{n}\longrightarrow\left(\Phi,0\right)$
converges regularly. On the other hand we have $t_{n}\left(0\right)=\dim_{H}\left(\Lambda_{\Phi^{n}}\right)$.
Hence, the following corollary is straightforward and may be viewed
as a generalisation of the continuity results in \cite{lambdatopologypreprint,MR2183304}
for the Hausdorff dimension of the limit sets.
\begin{cor}
[\bf {Continuity of Hausdorff dimension}] Let $\Phi^{n}=\left(\phi_{e}^{n}\right)_{e\in I_{n}}$,
$\Phi=\left(\phi_{e}\right)_{e\in I}$ be elements of $\mathrm{CIFS}\left(X\right)$
such that $\left(\Phi^{n},0\right)_{n}\longrightarrow\left(\Phi,0\right)$
converges regularly. Then \[
\lim_{n\to\infty}\dim_{H}\left(\Lambda_{\Phi^{n}}\right)=\dim_{H}\left(\Lambda_{\Phi}\right).\]

\end{cor}

\subsection*{Finite-to-infinite phase transition}

To complete the discussion of the Exhausting Principle we would like
to emphasise that the boundary values of the approximating spectra
in general do not converge to the corresponding value of the limiting
system, i.~e. we may have\begin{equation}
f_{n}\left(\alpha_{\pm}^{n}\right)\not\to f\left(\alpha_{\pm}\right).\label{eq:nonConverge}\end{equation}
We refer to the property of an infinite system having a discontinuity
of this kind in one of the boundary points as a \emph{finite-to-infinite
phase transition} \emph{in }$\alpha_{+}$, resp. $\alpha_{-}$. Let
us illustrate this property with the following concrete example. 
\begin{example}
[\bf {Gauss system}] \label{exa:GaussExp}Let $\Phi:=\left\{ \phi_{e}:x\mapsto1/\left(x+e\right):e\in\N\right\} $
denote the \emph{Gauss system} and the potential $\psi$ is given
by $\psi\left(\omega\right):=-2\log\omega_{1}$, for $\omega\in\N^{\N}$.
In \cite{JaerischKess:08} we have shown that the multifractal spectrum
is unimodal, defined on $\left[0,1\right]$, and in the boundary points
of the spectrum we have $f(0)=0$ and $f(1)=1/2$. Nevertheless, for
the exhausting systems $\Phi^{n}:=\left(\phi_{e}\right)_{1\leq e\leq n}$
and $\psi^{n}:=\psi\big|_{\left\{ 1,\ldots,n\right\} ^{\N}}$ we have
for their corresponding multifractal spectra $f_{n}\left(\alpha_{+}^{n}\right)=0$
for all $n\in\N$ giving rise to a finite-to-infinite phase transition
(see Fig \ref{fig:Exhausting-Principle}). A proof of this will be
postponed to the end of Section \ref{sec:Multifractal-analysis-for-cIFS}.
\begin{figure}[h]
\includegraphics[width=1\textwidth]{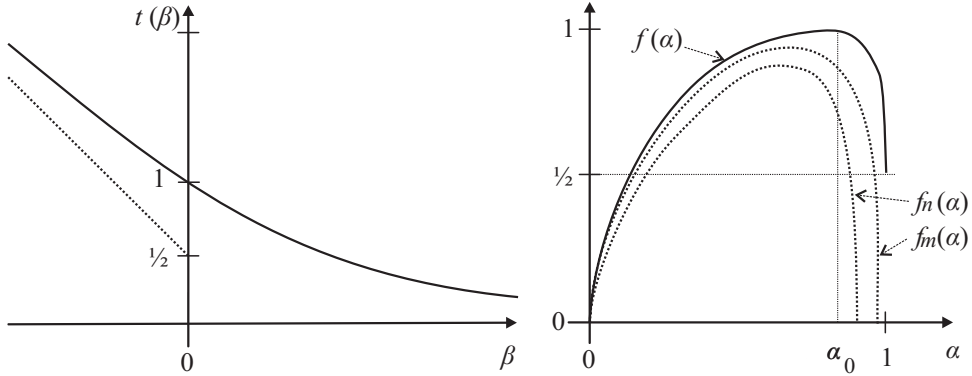}

\caption{\label{fig:Exhausting-Principle} Sketch illustrating the finite-to-infinite
phase transition for the Gauss system. The dashed graphs are associated
to the approximating spectra $f_{n}:\left[0,\alpha_{+}^{n}\right]\to\R_{+}$
and $f_{m}:\left[0,\alpha_{+}^{m}\right]\to\R_{+}$, $m>n$, of finite
sub-systems to the multifractal spectrum $f$ of the infinite system. }

\end{figure}

\end{example}
\begin{example}
[\bf {L\"uroth system}] \label{exa:-Lueroth} In the following example
the effective domain of the free energy function is not equal to $\R$,
which leads to an interesting boundary behaviour. For this let us
consider the \emph{L\"uroth system} $\Phi:=\left\{ \phi_{n}:x\mapsto x/(n(n+1))+1/\left(n+1\right):n\in\N\right\} $
(essentially a linearised Gauss system) and the potential functions
$\psi$ given by $\psi\left(\omega\right):=-\omega_{1}$, $\omega\in\N^{\N}$.
Then in virtue of our theorems the spectrum is given by the Legendre
transform of $t$ on $\left(2/\log\left(6\right),+\infty\right)$
via $f\left(\alpha\right)=-t^{*}\left(-\alpha\right)$. Similarly
as for the Gauss system in the example above, one can show that $f\left(2/\log\left(6\right)\right)=f_{n}\left(2/\log\left(6\right)\right)=f_{n}\left(\alpha_{+}^{n}\right)=0$,
$n\in\N$. Since we have Lebesgue almost everywhere that $\lim_{k\to\infty}\sum_{i=1}^{k}a_{i}/\sum_{i=1}^{k}\log\left(a_{i}\left(a_{i}+1\right)\right)=\int\left(\pi_{\Phi}^{-1}\left(x\right)\right)_{1}\,\mbox{d}\lambda/\int\log\left(\left(\pi_{\Phi}^{-1}\left(x\right)\right)_{1}\left(\left(\pi_{\Phi}^{-1}\left(x\right)\right)_{1}+1\right)\right)\,\mbox{d}\lambda=+\infty$
we find $f\left(+\infty\right)=1$. Hence as above, we have a finite-to-infinite
phase transition -- this time at infinity. 
\end{example}
\begin{example}
[\bf {Generalised L\"uroth system}] \label{exa:-LuerothGeneralized}
In the following example the effective domain of the free energy function
is again not equal to $\R$ and additionally we have a second order
phase transition. Let us consider the \emph{generalised} \emph{L\"uroth
system} $\Phi:=\left\{ \phi_{n}:x\mapsto4x/\left(n(n+1)(n+2)\right)+2/\left(\left(n+1\right)\left(n+2\right)\right):n\in\N\right\} $
and the potential functions $\psi$ given by $\psi\left(\omega\right):=-\omega_{1}$,
$\omega\in\N^{\N}$. Then in virtue of our theorems the spectrum is
given by the Legendre transform of $t$ on $\left(3/\log\left(15\right),\alpha_{+}\right)$
via $f\left(\alpha\right)=-t^{*}\left(-\alpha\right)$, where $\alpha_{+}:=\left(2\sum_{n\geq1}\frac{\log\left(\left(n(n+1)(n+2)\right)/4\right)}{\left(n(n+1)(n+2)\right)}\right)^{-1}$.
Using the Corollary \ref{cor:-ExhaustingPrinciple} (Exhausting Principle
II) (3) we gather some extra information on the spectrum. Since we
have $\lim_{n}-t_{n}^{*}\left(-t'_{n}\left(0\right)\right)=1$ and
$\alpha_{+}^{n}=n/\log(n(n+1)(n+2)/4)\to\infty$, we deduce that $1$
is a lower bound for $f\left(\alpha\right)$ for all $\alpha\geq\alpha_{+}$.
Similarly as for the Gauss system in the example above, one can show
that $f\left(3/\log\left(15\right)\right)=f_{n}\left(3/\log\left(15\right)\right)=f_{n}\left(\alpha_{+}^{n}\right)=0$,
$n\in\N$. Hence, $f\left(\alpha\right)=-t^{*}\left(-\alpha\right)$
for all $\alpha\geq3/\log(15)$. (cf. Fig. \ref{fig:Exhausting-Principle-ExpLueroth}).
\begin{figure}[h]
\includegraphics[width=1\textwidth]{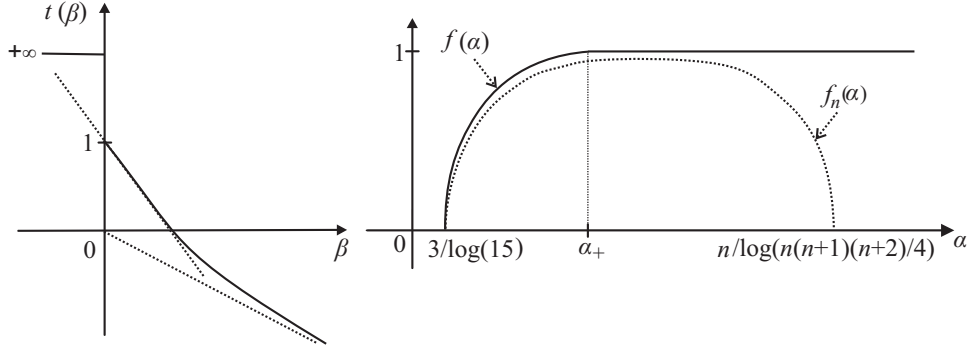}

\caption{\label{fig:Exhausting-Principle-ExpLueroth} Sketch of the free energy
function $t$ and the multifractal spectrum $f$ for the generalised
L\"uroth system. The dashed graph is associated to the approximating
spectra $f_{n}:\left[3/\log\left(15\right),n/\log(n(n+1)(n+2)/4)\right]\to\R_{+}$
of finite sub-system to the multifractal spectrum $f$ of the infinite
system. }

\end{figure}

Generalising further the latter two examples our analysis has successfully
been applied in \cite{KessMundayStatmann:10} to determine the Lyapunov
spectrum of $\alpha$-Farey-L\"uroth and $\alpha$-L\"uroth systems.
\end{example}
\begin{example}
[\bf {Irregular cIFS}] \label{exa:-IrregularIFS} For this example
we suppose that $\Phi$ is an irregular infinite cIFS, that is the
range of the pressure function $p:t\mapsto\mathcal{P}\left(t\zeta\right)$
consists of the negative reals and infinity (see \cite{urbanskimauldin-gdmsMR2003772}
for explicit examples), and let $\psi$ be constantly equal to $-1$.
We suppose that $p\left(\delta\right)=\eta<0$, where $\delta$ is
the critical value as well as the Hausdorff dimension of the limit
set. Then the free energy function $t$ is given by $t\left(\beta\right)=p^{-1}\left(\beta\right)$
for $\beta<\eta$ and constantly equal to $\delta$ for $\beta\geq\eta$.
The corresponding spectrum will have a linear part in $(0,\alpha_{-})$
if $-p^{+}\left(\delta\right)=1/\alpha_{-}<\infty$ and hence for
$\alpha_{-}>0$, we observe a \emph{second order phase transition}
(see Fig. \ref{fig:Exhausting-Principle-Critical}). %
\begin{figure}[h]
\includegraphics[width=1\textwidth]{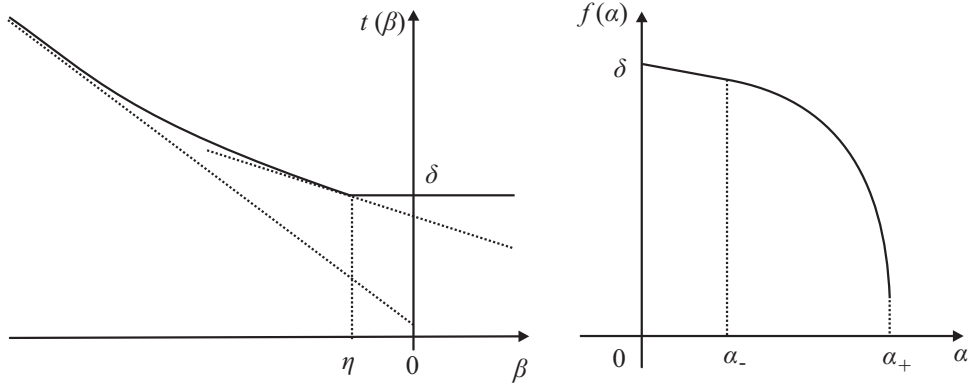}

\caption{\label{fig:Exhausting-Principle-Critical} Sketch of the free energy
function $t$ and the multifractal spectrum $f$ for an irregular
system with constant negative potential. Note that in this situation
we have a second order phase transition in $\alpha_{-}$ and the spectrum
$f$ is linear on $(0,\alpha_{-})$.}

\end{figure}

\end{example}
The paper is organised as follows. In Section \ref{sec:Preliminaries-from-CIFS}
we recall the basic notions relevant for cIFS. In Section \ref{sec:RegularityFreeEnergy}
we show the regularity of the free energy function proving Theorem
\ref{thm:convergence of temperature function}. Section \ref{sec:ConvergenceconjugacyConvexFunctions}
provides us with the necessary prerequisites from convex analysis
allowing us to deduce the multifractal regularity in Section \ref{sec:Multifractal-analysis-for-cIFS}.
In particular, we prove Theorem \ref{thm:tdach als hausdorff dim}
and \ref{thm:multifractal regularity}, and the finite-to-infinite
phase transition for the Gauss system. The final section is devoted
to the connection between our notion of regularity and the $\lambda$-topology.

\section{\textup{\label{sec:Preliminaries-from-CIFS} }Preliminaries }

Let us recall the definition of a conformal iterated function system
(see \cite{urbanskimauldin-gdmsMR2003772} for further details). Let
$X$ be a compact metric space. For an alphabet $I\subset\N$ with
$\card\left(I\right)\geq2$ we call $\Phi=\left(\phi_{e}\right)_{e\in I}$
an \emph{iterated function system} (IFS), where $\phi_{e}:X\rightarrow X$
are injective contractions, $e\in I$, with Lipschitz constants globally
bounded away from $1$.

Let $I^{*}:=\bigcup_{n\in\N}I^{n}$ denote the set of all finite subwords
of $I^{\N}$. We will consider the left shift map $\sigma:I^{\N}\rightarrow I^{\N}$
defined by $\sigma\left(\omega_{i}\right):=\left(\omega_{i+1}\right)_{i\geq1}$.
For $\omega\in I^{*}$ we let $|\omega|$ denote the length of the
word $\omega$, i.~e. the unique $n\in\N$ such that $\omega\in I^{n}$.

The space $I^{\N}$ is equipped with the metric $d$ given by $d(\omega,\tau):=\exp\left(-|\omega\wedge\tau|\right),$
where $\omega\wedge\tau\in I^{*}\cup I^{\N}$ denotes the longest
common initial block of the infinite words $\omega$ and $\tau$. 

We now describe the limit set of the iterated function system $\Phi$.
For each $\omega\in I^{*}$, say $\omega\in I^{n}$, we consider the
map coded by $\omega$,\[
\phi_{\omega}:=\phi_{\omega_{1}}\circ\cdots\circ\phi_{\omega_{n}}:X\rightarrow X.\]
For $\omega\in I^{\N}$, the sets $\left\{ \phi_{\omega|_{n}}\left(X\right)\right\} _{n\geq1}$
form a descending sequence of non-empty compact sets and therefore
$\bigcap_{n\geq1}\phi_{\omega|_{n}}\left(X\right)\neq\emptyset$.
Since for every $n\in\N$, $\diam\left(\phi_{\omega|_{n}}\left(X\right)\right)\leq s_{\Phi}^{n}\diam\left(X\right)$,
we conclude that the intersection $\bigcap\phi_{\omega|_{n}}\left(X\right)\in X$
is a singleton and we denote its only element by $\pi_{\Phi}(\omega)$.
In this way we have defined the coding map $\pi=\pi_{\Phi}:I^{\N}\rightarrow X.$
The set $\Lambda=\Lambda_{\Phi}=\pi\left(I^{\N}\right)$ will be called
the \emph{limit set} of $\Phi$.
\begin{defn}
\label{def:Conformal} We call an iterated function system \emph{conformal}
(cIFS) if the following conditions are satisfied. \renewcommand{\theenumi}{\alph{enumi}}
\begin{enumerate}
\item \label{enu:CIFSa}The \emph{phase space} $X$ is a compact connected
subset of a Euclidean space $\R^{D}$, $D\geq1$, such that $X$ is
equal to the closure of its interior, i.~e. $X=\overline{\Int(X)}$. 
\item \label{enu:CIFSb}\textit{\emph{(}}\textit{Open set condition }\textit{\emph{(OSC))}}
For all $a,b\in I$, $a\ne b$, \[
\phi_{a}\left(\Int(X)\right)\cap\phi_{b}\left(\Int(X)\right)=\emptyset.\]

\item \label{enu:CIFSc} There exists an open connected set $W\supset X$
such that for every $e\in I$ the map $\phi_{e}$ extends to a $C^{1}$
conformal diffeomorphism of $W$ into $W$. 
\item \label{enu:CIFSd}\textit{\emph{ (}}\textit{Cone property}\textit{\emph{)}}
There exist $\gamma,l>0$, $\gamma<\pi/2$, such that for every $x\in X\subset\R^{D}$
there exists an open cone $\Con(x,\gamma,l)\subset\Int(X)$ with vertex
$x$, central angle of measure $\gamma$, and altitude $l$. 
\item \label{enu:CIFSe} There are two constants $L=L_{\Phi}\geq1$ and
$\alpha=\alpha_{\Phi}>0$ such that \[
\left||\phi_{e}'(y)|-|\phi_{e}'(x)|\right|\leq\frac{L_{\Phi}}{\Vert(\phi_{e}')^{-1}\Vert_{X}}\Vert y-x\Vert^{\alpha}\]
 for every $e\in I$ and every pair of points $x,y\in X$.
\end{enumerate}
For the following let \begin{equation}
s_{\Phi}:=\sup_{e\in I}\left\Vert \left(\phi_{e}\right)'\right\Vert _{X}<1.\label{eq:def s_phi}\end{equation}
For a fixed phase space $X$ satisfying (\ref{enu:CIFSa}) the set
of conformal iterated function systems will be denoted \[
\mathrm{CIFS}\left(X\right):=\left\{ \Phi=\left(\phi_{e}:X\longrightarrow X\right)_{e\in I}\textrm{ cIFS, }I\subset\N\right\} .\]
 
\end{defn}
The following fact was proved in \cite{urbanskimauldin-gdmsMR2003772}.
\begin{prop}
\label{p1.033101} For $D\geq2$, any family $\Phi=\left(\phi_{e}\right)_{e\in I}$
satisfying condition \emph{(\ref{enu:CIFSa})} and \emph{(\ref{enu:CIFSc})}
also satisfies condition \emph{(\ref{enu:CIFSe})} with $\alpha=1$. 
\end{prop}
In \cite{urbanskimauldin-gdmsMR2003772} we also find the following
straightforward consequence of (\ref{enu:CIFSe}).
\begin{lem}
\label{l2.033101} If $\Phi=\left(\phi_{e}\right)_{e\in I}$ is a
cIFS, then for all $\omega\in I^{*}$ and all $x,y\in W$, we have
\[
\left|\log|\phi_{\omega}'(y)|-\log|\phi_{\omega}'(x)|\right|\leq\frac{L}{1-s^{\alpha}}\Vert y-x\Vert^{\alpha}.\]

\end{lem}
Another consequence of (\ref{enu:CIFSe}) is
\begin{itemize}
\item [(f)] \label{ite:cIFSf} (\emph{Bounded distortion property}). There
exists $K_{\Phi}\geq1$ such that for all $\omega\in I^{*}$ and all
$x,y\in X$ \[
|\phi_{\omega}'(y)|\leq K_{\Phi}|\phi_{\omega}'(x)|.\]

\end{itemize}
In \cite[Lemma 2.3.1]{urbanskimauldin-gdmsMR2003772} it has been
shown that for a H\"older continuous function $g:J^{\N}\rightarrow\R$,
$J\subset I$, we have for all $\omega\in J^{*}$ and all $x,y\in\left[\omega\right]$
that \[
\exp S_{\left|\omega\right|}g\left(x\right)\le K_{g}\exp S_{\left|\omega\right|}g\left(y\right).\]
In here, the constant $K_{g}\geq1$ only depends on the H\"older norm
and the H\"older exponent of $g$, as well as the metric on $J^{\N}$.
With \[
Z_{n}\left(g\right):=\sum_{\omega\in J^{n}}\exp\sup_{\tau\in[\omega]\cap J^{\N}}\left(S_{n}g\left(\tau\right)\right)\]
we will denote the \emph{$n$-th partition function} of $g$. 
\begin{defn}
\label{def:The-topological-pressure} The \emph{topological pressure}
$\P(f)$ of a continuous function $f:I^{\N}\rightarrow\R$ is defined
by the following limit, which always exists (possibly equal to $+\infty$),
\[
\P(f):=\lim_{n\rightarrow\infty}\frac{1}{n}\log Z_{n}\left(f\right)=\inf_{n}\frac{1}{n}\log Z_{n}\left(f\right).\]

\end{defn}
At the end of this section we would like to comment on the topology
of pointwise convergence. $\rho$ is well defined, since $\bigl\Vert\phi_{i}^{1}-\phi_{i}^{2}\bigr\Vert_{X}+\bigl\Vert(\phi_{i}^{1})'-(\phi_{i}^{2})'\bigr\Vert_{X}$
is bounded by $\diam\left(X\right)+2$. Using the fact that $A\triangle C\subset A\triangle B\cup B\triangle C$
for arbitrary sets $A,B,C$, we readily observe that $\rho$ as given
in (\ref{eq:RhoMetrik}) actually defines a metric on $\mathrm{CIFS}\left(X\right)$.
This metric induces the \emph{topology of pointwise convergence} on
$\mathrm{CIFS}\left(X\right)$. Let $\Phi^{n}=\left(\phi_{i}^{n}:X\longrightarrow X\right)_{i\in I_{n}},$
$n\in\N$, and $\Phi=\left(\phi_{i}:X\longrightarrow X\right)_{i\in I}$
be elements of $\mathrm{CIFS}\left(X\right)$ with $\Phi^{n}\rightarrow\Phi$
pointwise. Then for every $k\in\N$ we find an integer $N_{k}$ such
that and all $n\ge N_{k}$ we have $I_{n}\triangle I\subset\left\{ k+1,k+2,\dots\right\} .$
Similarly as in \cite[ Lemma 5.1]{MR2183304} it follows that pointwise
convergence in $\mathrm{CIFS}\left(X\right)$ is equivalent to the
following condition. 
\begin{condition}
\label{con:pointwiseconv-forgeometricpotential} We have $\left(\Phi^{n},\psi^{n}\right)\rightarrow\left(\Phi,\psi\right)$
pointwise if and only if for every $\omega\in I^{*}$ \[
\lim_{n\rightarrow\infty}\left(\bigl\Vert\phi_{\omega}^{n}-\phi_{\omega}\bigr\Vert_{X}+\bigl\Vert(\phi_{\omega}^{n})'-\phi_{\omega}'\bigr\Vert_{X}\right)=0\quad\mbox{and }\lim_{n\rightarrow\infty}\left(\bigl\Vert S_{|\omega|}\psi^{n}-S_{|\omega|}\psi\bigr\Vert_{\left[\omega\right]\cap I_{n}^{\N}}\right)=0\]

\end{condition}
Using this condition we are able to prove a technical property which
will be crucial in the proof of our main theorems. 
\begin{lem}
\label{lem:ConditionM*forZeta+Psi} Assume that $\left(\Phi^{n},\psi^{n}\right)\rightarrow\left(\Phi,\psi\right)$
converges pointwise. Then there exists $M>0$ such that for every
$\omega\in I^{*}$ fixed and all sufficiently large $n\in\N$ (depending
on $\omega$) we have for all $\eta\in\left[\omega\right]\cap I_{n}^{\N}$
and $\tau\in\left[\omega\right]\cap I^{\N}$ \[
\frac{\e^{S_{|\omega|}\psi^{n}\left(\eta\right)}}{\e^{S_{|\omega|}\psi\left(\tau\right)}},\frac{\e^{S_{|\omega|}\zeta^{n}\left(\eta\right)}}{\e^{S_{|\omega|}\zeta\left(\tau\right)}}\in\left[M^{-1},M\right].\]
\end{lem}
\begin{proof}
Using the above Condition \ref{con:pointwiseconv-forgeometricpotential}
we find for $\omega\in I^{*}$ and $n$ sufficiently large that\[
\max\left\{ \bigl\Vert\log\left|\left(\phi_{\omega}^{n}\right)'\right|-\log\left|\left(\phi_{\omega}\right)'\right|\bigr\Vert_{X},\bigl\Vert S_{|\omega|}\psi^{n}-S_{|\omega|}\psi\bigr\Vert_{\left[\omega\right]\cap I_{n}^{\N}}\right\} \le1.\]
Then we have with $K_{\Phi}$ and $K_{\psi}$ denoting the bounded
distortion constants as defined above \begin{align*}
\left|S_{|\omega|}\zeta^{n}\left(\eta\right)-S_{|\omega|}\zeta\left(\tau\right)\right| & =\left|\log\left|\left(\phi_{\omega}^{n}\right)'\right|\left(\pi_{\Phi^{n}}\left(\sigma^{n}\left(\eta\right)\right)\right)-\log\left|\left(\phi_{\omega}\right)'\right|\left(\pi_{\Phi}\left(\sigma^{n}\left(\tau\right)\right)\right)\right|\\
 & \le\left|\log\left|\left(\phi_{\omega}^{n}\right)'\right|\left(\pi_{\Phi^{n}}\left(\sigma^{n}\left(\eta\right)\right)\right)-\log\left|\left(\phi_{\omega}\right)'\right|\left(\pi_{\Phi^{n}}\left(\sigma^{n}\left(\eta\right)\right)\right)\right|\\
 & \qquad+\left|\log\left|\left(\phi_{\omega}\right)'\right|\left(\pi_{\Phi^{n}}\left(\sigma^{n}\left(\eta\right)\right)\right)-\log\left|\left(\phi_{\omega}\right)'\right|\left(\pi_{\Phi}\left(\sigma^{n}\left(\tau\right)\right)\right)\right|\\
 & \le1+\log K_{\Phi}\end{align*}
as well as \begin{align*}
\left|S_{|\omega|}\psi^{n}\left(\eta\right)-S_{|\omega|}\psi\left(\tau\right)\right| & \le\left|S_{|\omega|}\psi^{n}\left(\eta\right)-S_{|\omega|}\psi\left(\eta\right)\right|+\left|S_{|\omega|}\psi\left(\eta\right)-S_{|\omega|}\psi\left(\tau\right)\right|\\
 & \le1+\log K_{\psi}.\end{align*}
Letting $M:=3\cdot\max\left\{ K_{\Phi},K_{\psi}\right\} $ the lemma
follows. \end{proof}
\begin{rem}
Note that we may replace condition (\ref{enu:locallyUniformPsi})
in Definition \ref{def:pointwise convergence} by the slightly weaker
conditions on $\psi^{n}$ and $\psi$ stated in the above Lemma combined
with the condition that $\psi^{n}$ converges uniformly to $\psi$
on compact $\sigma$-invariant subsets of $I^{\N}$.
\end{rem}

\section{\label{sec:RegularityFreeEnergy}Regularity of the free energy function}

In this section we give a proof of Theorem \ref{thm:convergence of temperature function}.
Let $\zeta$ denote the geometric potential function associated with
$\Phi$ as defined in the introduction. For $\Phi\in\mathrm{CIFS}\left(X\right)$
and a H\"older continuous potential $\psi:I^{\N}\longrightarrow\R$
let $t$ denote the \emph{free energy function} of $\left(\Phi,\psi\right)$
as introduced in Definition \ref{def:temperature functions}, i.~e.
$t\left(\beta\right):=\inf\left\{ t:\P\left(t\zeta+\beta\psi\right)\leq0\right\} $.
Clearly, if there exists a zero of $t\longmapsto\P(t\zeta+\beta\psi)$
then $t\left(\beta\right)$ is the unique zero of this function (which
in particular is the case for a finite alphabet $I$). Also, $t\left(\beta\right)=+\infty$
if and only if $\left\{ t:\P\left(t\zeta+\beta\psi\right)\leq0\right\} =\emptyset$.
\begin{lem}
\label{lem:ProperConvex} The free energy $t$ of $\left(\Phi,\psi\right)$
is a proper (not necessarily closed) convex function on $\R$. \end{lem}
\begin{proof}
Fix $\beta_{1},\beta_{2}\in\R$, $\lambda\in\left(0,1\right)$ and
$\epsilon>0$. Using the convexity of the topological pressure we
have \begin{align*}
 & \P\left(\left(\lambda t\left(\beta_{1}\right)+\left(1-\lambda\right)t\left(\beta_{2}\right)+\epsilon\right)\zeta+\left(\lambda\beta_{1}+\left(1-\lambda\right)\beta_{2}\right)\psi\right)\\
 & \quad\quad=\P\left(\lambda\left(\left(t\left(\beta_{1}\right)+\epsilon\right)\zeta+\beta_{1}\psi\right)+\left(1-\lambda\right)\left(\left(t\left(\beta_{2}\right)+\epsilon\right)\zeta+\beta_{2}\psi\right)\right)\\
 & \quad\quad\le\lambda\P\left(\left(t\left(\beta_{1}\right)+\epsilon\right)\zeta+\beta_{1}\psi\right)+\left(1-\lambda\right)\P\left(\left(t\left(\beta_{2}\right)+\epsilon\right)\zeta+\beta_{2}\psi\right)\le0.\end{align*}
Hence, by definition of $t$, this implies $t\left(\lambda\beta_{1}+\left(1-\lambda\right)\beta_{2}\right)\le\lambda t\left(\beta_{1}\right)+\left(1-\lambda\right)t\left(\beta_{2}\right)+\epsilon$.
Since $\epsilon>0$ was arbitrary this shows the convexity. 

To see that $t$ is a proper convex function observe that $-\infty<\P\left(\beta\psi\right)$
and hence, for $t<0$, we have \[
\frac{1}{n}\log\sum_{\omega\in I^{n}}\exp\sup_{\tau\in\left[\omega\right]}\left(S_{n}t\zeta\left(\tau\right)+\beta\psi\right)\geq t\log\left(s_{\Phi}\right)+\P\left(\beta\psi\right)\to\infty\quad\mbox{ for }t\to-\infty.\]
Consequently, $t\left(\beta\right)>-\infty$ for all $\beta\in\R$.
Since also $t\left(0\right)=\dim_{H}\left(\Lambda_{\Phi}\right)<\infty$
we have that $t$ is proper. \end{proof}
\begin{lem}
\label{cla:ThemPressure} Let $\left(\Phi^{n},\psi^{n}\right)_{n}\longrightarrow\left(\Phi,\psi\right)$
converge regularly. Then for all $t,\beta\in\R$ and $n$ tending
to infinity we have \textup{\[
\P(t\zeta^{n}+\beta\psi^{n})\longrightarrow\P(t\zeta+\beta\psi).\]
} \end{lem}
\begin{proof}
Fix $t,\beta\in\R$ with $\P(t\zeta+\beta\psi)<\infty$ and $\epsilon>0$.
With $k\in\N$ and $C>0$ chosen according to Definition \ref{def:regularconvergence-1}
choose $m\in\N$ large enough such that $\left(mk\right)^{-1}\log Z_{mk}\left(t\zeta+\beta\psi\right)\le\P\left(t\zeta+\beta\psi\right)+\epsilon/2$
and $\left(mk\right)^{-1}\log\left(K+1\right)<\epsilon/2$, where
$K:=M^{\left|\beta\right|+|t|}$ and $M$ is the constant defined
in the proof of Lemma \ref{lem:ConditionM*forZeta+Psi}. We prove
that for all $n\in\N$ sufficiently large we have \begin{equation}
\left(mk\right)^{-1}\log Z_{mk}\left(t\zeta^{n}+\beta\psi^{n}\right)\le\left(mk\right)^{-1}\log Z_{mk}\left(t\zeta+\beta\psi\right)+\epsilon/2.\label{eq:uppersemi for Zk}\end{equation}
This would imply \begin{align*}
\P(t\zeta^{n}+\beta\psi^{n}) & \le\left(mk\right)^{-1}\log Z_{mk}\left(t\zeta^{n}+\beta\psi^{n}\right)\\
 & \le\left(mk\right)^{-1}\log Z_{mk}\left(t\zeta+\beta\psi\right)+\epsilon/2\le\P\left(t\zeta+\beta\psi\right)+\epsilon\end{align*}
for sufficiently large $n\in\N$. To prove (\ref{eq:uppersemi for Zk})
we first choose a finite set $F\subset I$ such that \[
\sum_{\omega\in I^{mk}\setminus F^{mk}}\exp\sup_{\rho\in\left[\omega\right]\cap I^{\N}}\left(S_{mk}\left(t\zeta+\beta\psi\right)\left(\rho\right)\right)<\frac{1}{\left(CK\right)^{m}}Z_{mk}\left(t\zeta+\beta\psi\right).\]
Then by Definition \ref{def:regularconvergence-1} and the choice
of $F$ we have \begin{align*}
 & \!\!\!\!\!\!\!\!\!\!\!\!\!\!\!\!\!\!\!\!\!\!\!\!\!\!\!\!\!\!\sum_{\omega\in I_{n}^{mk}\setminus F^{mk}}\exp\sup_{\rho\in I_{n}^{\N}\cap\left[\omega\right]}\left(S_{mk}\left(t\zeta^{n}+\beta\psi^{n}\right)\left(\rho\right)\right)\\
 & \leq\left(\sum_{\omega\in I_{n}^{k}\setminus F^{k}}\exp\sup_{\rho\in I_{n}^{\N}\cap\left[\omega\right]}\left(S_{k}\left(t\zeta^{n}+\beta\psi^{n}\right)\left(\rho\right)\right)\right)^{m}\\
 & \leq C^{m}\left(\sum_{\omega\in I_{n}^{k}\setminus F^{k}}\exp\sup_{\rho\in I^{\N}\cap\left[\omega\right]}\left(S_{k}\left(t\zeta+\beta\psi\right)\left(\rho\right)\right)\right)^{m}\\
 & \leq\left(CK\right)^{m}\sum_{\omega\in I_{n}^{mk}\setminus F^{mk}}\exp\sup_{\rho\in I^{\N}\cap\left[\omega\right]}\left(S_{mk}\left(t\zeta+\beta\psi\right)\left(\rho\right)\right)<Z_{mk}\left(t\zeta+\beta\psi\right)\end{align*}
Hence, on the one hand, we have \begin{align*}
Z_{mk}\left(t\zeta^{n}+\beta\psi^{n}\right) & =\sum_{\omega\in I_{n}^{mk}\cap F^{mk}}\exp\sup_{\rho\in I_{n}^{\N}\cap\left[\omega\right]}\left(S_{mk}\left(t\zeta^{n}+\beta\psi^{n}\right)\left(\rho\right)\right)\\
 & \quad\qquad\quad+\sum_{\omega\in I_{n}^{mk}\setminus F^{mk}}\exp\sup_{\rho\in I_{n}^{\N}\cap\left[\omega\right]}\left(S_{mk}\left(t\zeta^{n}+\beta\psi^{n}\right)\left(\rho\right)\right)\\
 & \le\sum_{\omega\in I_{n}^{mk}\cap F^{mk}}\exp\sup_{\rho\in I_{n}^{\N}\cap\left[\omega\right]}\left(S_{mk}\left(t\zeta^{n}+\beta\psi^{n}\right)\left(\rho\right)\right)+Z_{mk}\left(t\zeta+\beta\psi\right).\end{align*}
To find an upper bound also for the finite sum in the latter inequality
we note that by Lemma \ref{lem:ConditionM*forZeta+Psi} we have for
every $t\in\R$ and $\omega\in I_{n}^{mk}\cap F^{mk}$ and for sufficiently
large $n$ that \[
\exp\sup_{\tau\in I_{n}^{\N}\cap\left[\omega\right]}\left(S_{mk}\left(t\zeta^{n}+\beta\psi^{n}\right)\left(\tau\right)\right)\le K\exp\sup_{\rho\in I^{\N}\cap\left[\omega\right]}\left(S_{mk}\left(t\zeta+\beta\psi\right)\left(\rho\right)\right).\]
Since $I_{n}^{mk}\cap F^{mk}$ is finite we have on the other hand
for $n$ sufficiently large that \begin{align*}
\sum_{\omega\in I_{n}^{mk}\cap F^{mk}}\e^{\sup_{\rho\in I_{n}^{\N}\cap\left[\omega\right]}\left(S_{mk}\left(t\zeta^{n}+\beta\psi^{n}\right)\left(\rho\right)\right)} & \le K\sum_{\omega\in I_{n}^{mk}\cap F^{mk}}\e^{\sup_{\rho\in I^{\N}\cap\left[\omega\right]}\left(S_{mk}\left(t\zeta+\beta\psi\right)\left(\rho\right)\right)}.\end{align*}
Combining both estimates we find for $n$ sufficiently large \begin{eqnarray*}
Z_{mk}\left(t\zeta^{n}+\beta\psi^{n}\right) & \le & \left(K+1\right)Z_{mk}\left(t\zeta+\beta\psi\right).\end{eqnarray*}
Taking logarithm and dividing by $mk$ proves (\ref{eq:uppersemi for Zk}). 

To prove the reverse inequality $\liminf_{n}\P(t\zeta^{n}+\beta\psi^{n})\ge\P(t\zeta+\beta\psi)$
let us fix $\epsilon>0$. Using \cite[Theorem 2.15]{urbanskimauldin-gdmsMR2003772}
we can choose a finite set $F\subset I$ such that \[
\P(t\zeta\big|_{F^{\N}}+\beta\psi\big|_{F^{\N}})\ge\P(t\zeta+\beta\psi)-\epsilon.\]
By \cite[Lemma 4.2]{MR2183304} and Definition \ref{def:pointwise convergence}
(\ref{enu:locallyUniformPsi}) we have $\bigl\Vert\zeta^{n}\big|_{F^{\N}}-\zeta\big|_{F^{\N}}\bigr\Vert_{F^{\N}}\rightarrow0$
as well as $\bigl\Vert\psi^{n}\big|_{F^{\N}}-\psi\big|_{F^{\N}}\bigr\Vert_{F^{\N}}\rightarrow0$.
Since $f\mapsto\P\left(f\right)$ is Lipschitz-continuous with respect
to the $\Vert\cdot\Vert_{F^{\N}}$-norm (cf. \cite[Theorem 9.7]{walters-ergodictheoryMR648108})
we conclude \begin{eqnarray*}
\liminf_{n}\P(t\zeta^{n}+\beta\psi^{n}) & \ge & \liminf_{n}\P(t\zeta^{n}\big|_{F^{\N}}+\beta\psi^{n}\big|_{F^{\N}})\\
 & = & \P(t\zeta\big|_{F^{\N}}+\beta\psi\big|_{F^{\N}})\\
 & \ge & \P(t\zeta+\beta\psi)-\epsilon.\end{eqnarray*}

\end{proof}
\begin{proof}
[Proof of\emph{ }Theorem \ref{thm:convergence of temperature function}]
Fix $\beta\in\R$. To verify $\limsup_{n}t_{n}\left(\beta\right)\le t\left(\beta\right)$
we may assume $t\left(\beta\right)<\infty$. Since the map $p_{\beta}:t\mapsto\P(t\zeta+\beta\psi)$
is strictly decreasing on $\dom p_{\beta}$ we have that $\P\left(\left(t\left(\beta\right)+\delta\right)\zeta+\beta\psi\right)<0$
for every $\delta>0$. As a consequence of Lemma \ref{cla:ThemPressure}
we have $\P\left(\left(t\left(\beta\right)+\delta\right)\zeta^{n}+\beta\psi^{n}\right)<0$
for all $n$ sufficiently large. This implies $t_{n}\left(\beta\right)\le t\left(\beta\right)+\delta$
for all $\delta>0$ and therefore $\limsup_{n}t_{n}\left(\beta\right)<t\left(\beta\right)$. 

To verify $\liminf_{n}t_{n}\left(\beta\right)\ge t\left(\beta\right)$
we first assume that $t\left(\beta\right)<\infty.$ By definition
of $t$ we have $\P\left(\left(t\left(\beta\right)-\delta\right)\zeta+\beta\psi\right)>0$
for every $\delta>0$. Then again by Lemma \ref{cla:ThemPressure}
we also have $\P\left(\left(t\left(\beta\right)-\delta\right)\zeta^{n}+\beta\psi^{n}\right)>0$
for all $n$ sufficiently large, which in turn implies $t_{n}\left(\beta\right)\ge t\left(\beta\right)-\delta$
for all $n$ large enough. Finally, let $t\left(\beta\right)=\infty$,
i.~e. $\P\left(t\zeta+\beta\psi\right)=\infty$ for all $t$. By
Lemma \ref{cla:ThemPressure} we have for any $t\in\R$ fixed $\P\left(t\zeta^{n}+\beta\psi^{n}\right)>0$
for $n$ large enough and hence $t_{n}\left(\beta\right)\ge t$ for
$n$ large enough. Since $t\in\R$ was arbitrary it follows that $t_{n}\left(\beta\right)$
tends to infinity as $n$ increases.
\end{proof}

\section{\label{sec:ConvergenceconjugacyConvexFunctions}Convergence and conjugacy
of convex functions}

In this section we collect the necessary basic facts from convex analysis
needed for the multifractal analysis in Section \ref{sec:Multifractal-analysis-for-cIFS}.
We closely follow \cite{MR0479398}, and all details can be found
either therein or in \cite{rockafellar-convexanalysisMR0274683}.

The following proposition is a direct consequence of \cite[Corollaries 2C and 3B]{MR0479398}
combined with the fact that Legendre conjugation is continuous with
respect to the convergence of epigraphs in the classical sense as
defined e.~g. by Kuratowski in \cite{Kuratowski:66}. 
\begin{prop}
\label{pro:sawe pointwise convergence of legendre} Let $g_{n},g$,
$n\in\N$, be closed convex functions on $\R$ such that $\Int\left(\dom\left(g\right)\right)\neq\emptyset$
and $g_{n}\longrightarrow g$ pointwise. Then pointwise on $\Int\left(\dom\left(g^{*}\right)\right)$,
we have $g_{n}^{*}\longrightarrow g^{*}$. 
\end{prop}
The following corollary allows us to apply Proposition \ref{pro:sawe pointwise convergence of legendre}
also in the case when the functions $g_{n},g$ are not closed.
\begin{cor}
\label{cor:convergence conjugate functions}Let $g_{n},g$, $n\in\N$,
be convex functions on $\R$ and $\alpha\in\R$ such that there exist
$x_{1},x_{2}\in\Int\left(\dom\left(g\right)\right)$ with $x_{1}<x_{2}$
and $g^{+}\left(x_{1}\right)<\alpha<g^{-}\left(x_{2}\right)$. Furthermore,
assume that there exists an open neighbourhood $U\subset\dom\left(g\right)$
containing $x_{1},x_{2}$ such that $g_{n}\big|_{U}\longrightarrow g\big|_{U}$
pointwise. Then we have $g_{n}^{*}\left(\alpha\right)\longrightarrow g^{*}\left(\alpha\right)$.\end{cor}
\begin{proof}
Without loss of generality we have $g_{n}\big|_{U}<\infty$ for all
$n\in\N$. Let $\1_{A}$ denote the indicator function on the set
$A$ and let $\tilde{g}_{n},\tilde{g}$ denote the closed convex functions
given by $\tilde{g}_{n}:=g_{n}\1_{U}+\infty\1_{\R\setminus U}$ and
$\tilde{g}:=g\1_{U}+\infty\1_{\R\setminus U}$. Notice that these
closed convex functions agree on $U$ with the original functions.
By Proposition \ref{pro:sawe pointwise convergence of legendre} we
conclude that $\tilde{g}_{n}^{*}\longrightarrow\tilde{g}^{*}$ pointwise
on $\Int\left(\dom\tilde{g}^{*}\right)$. Clearly by our assumptions,
$\alpha\in\Int\left(\dom\tilde{g}^{*}\right)$ and $\alpha$ belongs
to the subdifferential $\partial g\left(x\right):=\left\{ a\in\R:\forall x'\in\R\; g(x')-g(x)\geq a\left(x'-x\right)\right\} $
for some $x\in U$, and hence $\tilde{g}^{*}\left(\alpha\right)=\alpha x-\tilde{g}\left(x\right)=g^{*}\left(\alpha\right)$
by \cite[Theorem 23.5]{rockafellar-convexanalysisMR0274683}. It remains
to show that $\tilde{g}_{n}^{*}\left(\alpha\right)=g_{n}^{*}\left(\alpha\right)$
for $n$ sufficiently large. Since by \cite[Theorem 24.5]{rockafellar-convexanalysisMR0274683}
the subdifferentials converge, the assumption $g^{+}\left(x_{1}\right)<\alpha<g^{-}\left(x_{2}\right)$
implies that $\alpha\in\partial g_{n}\left(y_{n}\right)$ for some
$y_{n}\in U$ and $n$ large. Then again by \cite[Theorem 23.5]{rockafellar-convexanalysisMR0274683}
we have $\tilde{g}_{n}^{*}\left(\alpha\right)=\alpha y_{n}-g_{n}\left(y_{n}\right)=g_{n}^{*}\left(\alpha\right)$. 
\end{proof}

\section{Regularity of the multifractal spectrum\label{sec:Multifractal-analysis-for-cIFS}}

We proceed by proving the Theorems \ref{thm:tdach als hausdorff dim}
and \ref{thm:multifractal regularity}. Recall that throughout we
use the generalised version of the free energy function $t$ as stated
in Definition \ref{def:temperature functions}. 
\begin{proof}
[Proof of Theorem \ref{thm:tdach als hausdorff dim}] Using the
definition of topological pressure and a standard covering argument
(just cover $\mathcal{F}_{\alpha}$ with cylinder sets) we obtain
$f\left(\alpha\right)\le\max\left\{ -t^{*}\left(-\alpha\right),0\right\} $
for every $\alpha\in\R$. 

We will use the \emph{Exhausting Principle} to prove the reverse inequality.
Let $\Phi^{n}=\left(\phi_{e}\right)_{e\in I_{n}}$ with $I_{n}:=I\cap\left\{ 1,\dots n\right\} $
and $\psi^{n}:=\psi\big|_{I_{n}^{\N}}$, $n\in\N$. Clearly, $\left(\Phi^{n},\psi^{n}\right)_{n}\longrightarrow\left(\Phi,\psi\right)$
converges regularly (Example \ref{exa:ExhaustionPrinciple}) and hence
by Theorem \ref{thm:convergence of temperature function}, we conclude
that $t_{n}\longrightarrow t$ pointwise on $\R$. Note that for $\alpha\in\left(\alpha_{-},\alpha_{+}\right)$
we find $x_{1},x_{2}\in\Int\left(\dom\left(t\right)\right)$ with
$x_{1}>x_{2}$ and $-t^{-}\left(x_{1}\right)<\alpha<-t^{+}\left(x_{2}\right)$.
Hence by Corollary \ref{cor:convergence conjugate functions}, we
conclude \[
\lim_{n\to\infty}-t_{n}^{*}\left(-\alpha\right)=-t^{*}\left(-\alpha\right).\]
Since the functions $t_{n}$ are finite and differentiable on $\R$
we conclude by \cite[Theorem 24.5]{rockafellar-convexanalysisMR0274683}
that $\alpha\in-t_{n}'\left(\R\right)$ for all $n$ large enough.
Recall that in the finite alphabet case it is well-known that $f_{n}\left(\alpha\right)=-t_{n}^{*}\left(-\alpha\right)$.
By construction we have $\mathcal{F}_{\alpha}^{n}\subset\mathcal{F}_{\alpha}$
and hence \[
-t^{*}\left(-\alpha\right)\geq f\left(\alpha\right)\geq f_{n}\left(\alpha\right)=-t_{n}^{*}\left(-\alpha\right)\;\longrightarrow\;-t^{*}\left(-\alpha\right).\]
 
\end{proof}
\begin{proof}
[Proof of Theorem \ref{thm:multifractal regularity}] By Theorem
\ref{thm:convergence of temperature function} the free energy functions
$t_{n}$ converge pointwise to $t$ on $\R$. For $\alpha\in\left(\alpha_{-},\alpha_{+}\right)$
we have by Corollary \ref{cor:convergence conjugate functions} that
$\lim_{n}-t_{n}^{*}\left(-\alpha\right)=-t^{*}\left(-\alpha\right)$.
Furthermore by Theorem \ref{thm:tdach als hausdorff dim}, we have
$-t^{*}\left(-\alpha\right)=f\left(\alpha\right)$. We also have $-t_{n}^{*}\left(-\alpha\right)=f_{n}\left(\alpha\right)$
for large $n$ by Theorem \ref{thm:tdach als hausdorff dim}, since
\cite[Theorem 24.5]{rockafellar-convexanalysisMR0274683} implies
$\alpha\in\left(\alpha_{-}^{n},\alpha_{+}^{n}\right)$ for large $n$.
This proves the first part of the theorem.

We now consider the case $\sup\dom\left(t\right)=+\infty$ (the second
case is proved along the same lines). Let us assume on the contrary
that there exists $\epsilon>0$ with $\alpha_{-}^{n}<\alpha_{-}-\epsilon$
for infinitely many $n$. For fixed $K>0$ we find $x\in\dom\left(t\right)$
with $x>\left(K+1\right)/\epsilon$ and such that $t'\left(x\right)$
exists. By Theorem \ref{thm:convergence of temperature function}
we find $n\in\N$ such that $\left\vert t_{n}\left(x\right)-t\left(x\right)\right\vert <1$
and $\alpha_{-}^{n}<\alpha_{-}-\epsilon$. Furthermore, we can choose
$\alpha\in\left(\alpha_{-}^{n},\alpha_{+}^{n}\right)$ satisfying
$\alpha<\alpha_{-}-\epsilon$. By Theorem \ref{thm:tdach als hausdorff dim}
we have $f_{n}\left(\alpha\right)=-t_{n}^{*}\left(-\alpha\right)$.
Since $\alpha+t'\left(x\right)\leq\alpha-\alpha_{-}<-\epsilon$ we
have by \cite[Theorem 23.5]{rockafellar-convexanalysisMR0274683}
\begin{align*}
f_{n}\left(\alpha\right)+t^{*}\left(t'\left(x\right)\right) & =-t_{n}^{*}\left(-\alpha\right)+t^{*}\left(t'\left(x\right)\right)\\
 & =\inf_{\beta}\left\{ t_{n}\left(\beta\right)+\beta\alpha\right\} -t\left(x\right)+xt'\left(x\right)\\
 & \leq\left(t_{n}\left(x\right)+x\alpha\right)-t\left(x\right)+xt'\left(x\right)\\
 & \leq1+x\left(\alpha+t'\left(x\right)\right)<1-\epsilon x<-K.\end{align*}
 Since $t^{*}\left(t'\left(x\right)\right)\ge-t\left(0\right)=-\dim_{H}\left(\Lambda\right)$
and $K$ can be chosen arbitrary large we get a contradiction to $f_{n}\left(\alpha\right)\geq0$. \end{proof}
\begin{rem}
\label{rem:ExhaustionPrinciple} Assume in the situation of Example
\ref{exa:ExhaustionPrinciple} that $\left(\alpha_{-},\alpha_{+}\right)$
is a bounded interval. If $t_{n}^{*}\left(\alpha_{\pm}\right)\longrightarrow t^{*}\left(\alpha_{\pm}\right)$
then $f_{n}\left(\alpha_{\pm}\right)\longrightarrow f\left(\alpha_{\pm}\right)$
and we have $-t^{*}\left(-\alpha_{\pm}\right)=f\left(\alpha_{\pm}\right)$.
To see this notice that $a\mapsto-t^{*}\left(-a\right)$ is bounded
on $\left(\alpha_{-},\alpha_{+}\right)$ since it coincides with the
Hausdorff dimension of certain sets. Furthermore, by \cite[Theorem 12.2]{rockafellar-convexanalysisMR0274683}
we have that $t^{*}$ is closed and hence $-t^{*}\left(-\alpha_{\pm}\right)<\infty$.
By our assumption we have $-t_{n}^{*}\left(-\alpha_{\pm}\right)<\infty$
for $n$ large, hence $\alpha_{\pm}\in\overline{-t_{n}'\left(\R\right)}$.
Then the claim follows by observing that \[
-t_{n}^{*}\left(-\alpha_{\pm}\right)=f_{n}\left(\alpha_{\pm}\right)\le f\left(\alpha_{\pm}\right)\le-t^{*}\left(-\alpha_{\pm}\right).\]

\end{rem}
Finally, we will sketch the proof of (\ref{eq:nonConverge}) for the
Gauss system as announced in the introduction. Let $\underline{n}:=\left(n,n,n,\ldots\right)$,
$n\geq2$, and suppose that $\omega\in\left\{ 1,\dots,n\right\} ^{k}$,
$k\in\N$, differs from $\underline{n}\big|_{k}$ in at least $\ell$
positions. Let $q_{k}\left(\mbox{\ensuremath{\omega}}\right)$ denote
the denominator of the $k$'s approximant of the continued fraction
expansion $[\omega_{1},\ldots,\omega_{k}]\in[0,1]$. Then, by the
recursive definition of $q_{k}$, we have \begin{equation}
\frac{q_{k}\left(\omega\right)}{q_{k}\left(\underline{n}\right)}>\left(\frac{n-1}{n}\right)^{\ell}.\label{eq:q_k}\end{equation}
 From this it follows that \[
\alpha_{+}^{n}=\frac{\psi\left(\underline{n}\right)}{\zeta\left(\underline{n}\right)}=\lim_{k\to\infty}\frac{k\log\left(n\right)}{\log q_{k}\left(\underline{n}\right)}=\frac{-\lg\left(n\right)}{\lg\left(-n/2+\sqrt{n^{2}/4+1}\right)}<1\]
 (cf. \cite[Fact 3]{JaerischKess:08}). Now we will argue similar
as in \cite{kessemultifractalsternbrocotMR2338129}. To prove that
$f_{n}\left(\alpha_{+}^{n}\right)=0$ it is sufficient to verify that
\[
\left\{ \mu\in\mathcal{M}\big(\left\{ 1,\ldots,n\right\} ^{\N},\sigma\big):\int\psi\,\mbox{d}\mu/\int\zeta\,\mbox{d}\mu=\alpha_{+}^{n}\right\} =\left\{ \delta_{\underline{n}}\right\} ,\]
where $\mathcal{M}\big(\left\{ 1,\ldots,n\right\} ^{\N},\sigma\big)$
denotes the set of shift invariant measures and $\delta_{x}$ the
Dirac measure centred on $x$. For the detailed argument see \cite{kessemultifractalsternbrocotMR2338129}.
We prove this fact by way of contradiction. Assume there exists $\mu\in\mathcal{M}\left(\left\{ 1,\ldots,n\right\} ^{\N},\sigma\right)$
with $\mu\not=\delta_{\underline{n}}$. By convexity of the set of
measures under consideration we may assume that $\mu$ is ergodic.
Then by our assumption there exists $\ell<n$ with $\mu\left([\ell]\right)=\eta>0$.
Then for all $\mu$-typical points $\omega$ we have $\lim_{k}S_{k}\psi\left(\omega\right)/S_{k}\zeta\left(\omega\right)=\lim_{k}\sum_{i=1}^{k}\log\omega_{i}/\log q_{k}\left(\omega\right)=\alpha_{+}^{n}$
and $S_{k}\1_{[\ell]}\left(\omega\right)\geq\eta k/2$ for all $k\in\N$
sufficiently large. Hence, using (\ref{eq:q_k}), we obtain\begin{align*}
 & \negmedspace\negmedspace\negmedspace\negmedspace\negmedspace\negmedspace\negmedspace\negmedspace\negmedspace\negmedspace\negmedspace\frac{k\log\left(n\right)}{\log q_{k}\left(\underline{n}\right)}-\frac{\sum^{k}\log\omega_{i}}{\log q_{k}\left(\omega\right)}=\\
 & =\frac{k\log\left(n\right)\log q_{k}\left(\omega\right)-\log q_{k}\left(\underline{n}\right)\sum^{k}\log\omega_{i}}{\log q_{k}\left(\omega\right)\log q_{k}\left(\underline{n}\right)}\\
 & =\frac{\left(\sum_{i=1}^{k}\log\left(n\right)-\log\left(\omega_{i}\right)\right)\log q_{k}\left(\omega\right)+\left(\log q_{k}\left(\omega\right)-\log q_{k}\left(\underline{n}\right)\right)\sum^{k}\log\omega_{i}}{\log q_{k}\left(\omega\right)\log q_{k}\left(\underline{n}\right)}\\
 & \geq\frac{k\left(\eta/2\right)\left(\log\left(n\right)-\log\left(n-1\right)\right)\log q_{k}\left(\omega\right)+\left(\log q_{k}\left(\omega\right)-\log q_{k}\left(\underline{n}\right)\right)\sum_{i=1}^{k}\log\omega_{i}}{\log q_{k}\left(\omega\right)\log q_{k}\left(\underline{n}\right)}\\
 & \geq\frac{k\eta\left(\log\left(n\right)-\log\left(n-1\right)\right)\left(\log q_{k}\left(\omega\right)-\sum_{i=1}^{k}\log\omega_{i}\right)}{2\log q_{k}\left(\omega\right)\log q_{k}\left(\underline{n}\right)}\\
 & \to\frac{\eta\alpha_{+}^{n}}{2\log n}\left(\log\left(n\right)-\log\left(n-1\right)\right)\left(1-\alpha_{+}^{n}\right)>0.\end{align*}
Since $\lim_{k\to\infty}k\log\left(n\right)/\log q_{k}\left(\underline{n}\right)=\alpha_{+}^{n}$
we obtain a contradiction.

\section{The extended $\lambda$-topology \label{sec:Connection with lambda}}

In this section we compare the notion of regular convergence with
the $\lambda$-topology introduced by Roy and Urba\'{n}ski. In particular,
as a consequence of Proposition \ref{pro:connection with lambda topology}
we will verify Example \ref{exa:lambda topology}. 

For ease of notation we will always assume $I=\N$. Let us first recall
the definition of the $\lambda$-topology from \cite{MR2183304} and
then give a generalisation to adapt this concept to our purposes.

For $\Phi^{n}=\left(\phi_{e}^{n}\right)_{e\in I}$, $\Phi=\left(\phi_{e}\right)_{e\in I}$
elements of $\mathrm{CIFS}\left(X\right)$ sharing the same alphabet
$I$ we say that $\Phi^{n}$ \emph{converges to $\Phi$ in the $\lambda$-topology,}
if $\Phi^{n}\rightarrow\Phi$ in the $\rho$-metric and there exists
$R>1$ such that for all sufficiently large $n$ and all $e\in I$
we have \begin{equation}
R^{-1}\le\frac{\bigl\Vert(\phi_{e}^{n})'\bigr\Vert_{X}}{\bigl\Vert\phi_{e}'\bigr\Vert_{X}}\le R.\label{eq:lambda topology property}\end{equation}
We shall generalise this to the case where we have $I_{n}\subset I=\N$.
We say that $\Phi^{n}=\left(\phi_{e}^{n}\right)_{e\in I_{n}}$ \emph{converges
to $\Phi=\left(\phi_{e}\right)_{e\in I}$ in the extended $\lambda$-topology}
of $\mathrm{CIFS}\left(X\right)$, if they converge in the $\rho$-metric
and there exists $D>1$ such that for all $n$ sufficiently large
and all $e\in I_{n}$ the assumption (\ref{eq:lambda topology property})
holds. 

Let us begin with the following basic lemma.
\begin{lem}
\label{lem:s_n converges to s }Let $\Phi^{n}=\left(\phi_{e}^{n}\right)_{e\in I_{n}}$,
$\Phi=\left(\phi_{e}\right)_{e\in\N}$ be elements of $\mathrm{CIFS}\left(X\right)$
with $\Phi^{n}\rightarrow\Phi$ converging in the \emph{extended}
$\lambda$-topology. Then with $s_{\Phi}$ defined in \emph{(}\textup{\ref{eq:def s_phi}}\emph{),
we have} \[
\lim_{n}s_{\Phi^{n}}=s_{\Phi}<1.\]
\end{lem}
\begin{proof}
By the open set condition (OSC) there exists $e\in I$ with $\bigl\Vert\phi_{e}'\bigr\Vert_{X}=\sup_{e}\bigl\Vert\phi_{e}'\bigr\Vert_{X}=s_{\Phi}$
as well as for every $n\in\N$ there exists $e_{n}\in I_{n}$ satisfying
$\bigl\Vert(\phi_{e_{n}}^{n})'\bigr\Vert_{X}=\sup_{e}\bigl\Vert(\phi_{e}^{n})'\bigr\Vert_{X}=s_{\Phi^{n}}$.
Since $\lim_{n}\bigl\Vert(\phi_{e}^{n})'\bigr\Vert_{X}=\bigl\Vert\phi_{e}'\bigr\Vert_{X}$
we have $\liminf_{n}s_{\Phi^{n}}\ge s_{\Phi}$. 

Next we conclude that $\left\{ e_{n}:n\in\N\right\} $ is contained
in a finite set $F\subset\N$. This follows by way of contradiction.
Assume the set is infinite. Then there exists a subsequence $n_{k}$
such that on the one hand $\bigl\Vert(\phi_{e_{n_{k}}}')\bigr\Vert_{X}\rightarrow0$
and on the other hand $\bigl\Vert(\phi_{e_{n_{k}}}^{n_{k}})'\bigr\Vert_{X}\rightarrow\liminf_{n}s_{\Phi^{n}}\geq s_{\Phi}>0$.
This would contradict property (\ref{eq:lambda topology property})
defining the extended $\lambda$-topology. Now by the definition of
the $\rho$-metric we have for all $\ell\in F$ that $\lim_{n}\bigl\Vert(\phi_{\ell}^{n})'\bigr\Vert_{X}=\bigl\Vert(\phi_{\ell})'\bigr\Vert_{X}\le s_{\Phi}$.
This gives $\limsup_{n}s_{\Phi^{n}}\le s_{\Phi}$. 
\end{proof}
For the following let $\psi^{n}:I_{n}^{\N}\longrightarrow\R$ and
$\psi:I^{\N}\longrightarrow\R$ be H\"older continuous functions, satisfying
condition (\ref{enu:locallyUniformPsi}) in Definition \ref{def:pointwise convergence}. 
\begin{assumption}
\label{ass:AdditionalBoundOnPsi }Additionally, we assume that there
exist $k\in\N$ and $M\in\N$, such that for all $n\in\N$ and for
all $\omega\in I_{n}^{k}$, $\tau\in I_{n}^{\N}\cap\left[\omega\right]$
and $\eta\in I^{\N}\cap\left[\omega\right]$ we have \[
M^{-1}\le\frac{\exp\left(S_{k}\left(\psi^{n}\right)\left(\tau\right)\right)}{\exp\left(S_{k}\left(\psi\right)\left(\eta\right)\right)}\le M.\]
\end{assumption}
\begin{rem}
Assumption \ref{ass:AdditionalBoundOnPsi } is for instance satisfied,
if 
\begin{enumerate}
\item $\sup_{n}\left\Vert \psi^{n}\right\Vert _{I_{n}^{\N}}<\infty$ and
$\left\Vert \psi\right\Vert _{I^{\N}}<\infty$, or 
\item $\sup_{n}\left\Vert \psi^{n}-\psi\right\Vert _{I_{n}^{\N}}<\infty$. 
\end{enumerate}
\end{rem}
For the following proposition recall that $K_{\Phi}$ denotes the
bounded distortion constant for $\Phi$ as stated in condition (f)
of the definition of a cIFS.
\begin{prop}
\label{pro:connection with lambda topology}Let $\Phi^{n}=\left(\phi_{e}^{n}\right)_{e\in I_{n}}$,
$\Phi=\left(\phi_{e}\right)_{e\in\N}$ be elements of $\mathrm{CIFS}\left(X\right)$
with $\Phi^{n}\rightarrow\Phi$ converging in the extended $\lambda$-topology.
Let $\psi^{n}:I_{n}^{\N}\longrightarrow\R$ and $\psi:I^{\N}\longrightarrow\R$
be H\"older continuous functions satisfying condition (\ref{enu:locallyUniformPsi})
in Definition \ref{def:pointwise convergence} as well as Assumption
\ref{ass:AdditionalBoundOnPsi }. Then $\left(\Phi^{n},\psi^{n}\right)_{n}\rightarrow\left(\Phi,\psi\right)$
converges regularly, if one of the following conditions is satisfied:
\begin{enumerate}
\item \label{enu:condition2}$\sup_{n}K_{\Phi^{n}}<\infty$. 
\item \label{enu:condition1}$\left\Vert \psi\right\Vert _{I^{\N}}<\infty$.
\item \label{enu:condition3} $\inf_{n}\alpha_{\Phi_{n}}>0$ and $\sup_{n}L_{\Phi_{n}}<\infty$.
\item \label{enu:condition4} $D\ge2$ and the maps $\Phi_{e}^{n}$ extend
to conformal diffeomorphisms on a common neighbourhood $W\supset X$
into $W$ for all $e\in I_{n}$ and $n\in\N$.
\end{enumerate}
\end{prop}
\begin{proof}
Clearly $\left(\Phi^{n},\psi^{n}\right)_{n}\rightarrow\left(\Phi,\psi\right)$
converges pointwise. Hence, we are left to verify the condition in
Definition \ref{def:regularconvergence-1} under the assumption (\ref{enu:condition2})
as well as under the assumption (\ref{enu:condition1}), and then
show that both (\ref{enu:condition3}) and (\ref{enu:condition4})
imply (\ref{enu:condition2}). 

\textbf{ad (\ref{enu:condition2}):} For $t\ge0$ we argue as follows.
Let $k\in\N$, $\omega\in I_{n}^{k}$ and $n$ sufficiently large,
such that (\ref{eq:lambda topology property}) and Assumption \ref{ass:AdditionalBoundOnPsi }
hold. Using this and the bounded distortion property of $\Phi$ from
Lemma \ref{l2.033101} with bounded distortion constant $K=K_{\Phi}$
we obtain\begin{align*}
\exp\sup_{\tau\in I_{n}^{\N}\cap\left[\omega\right]}\left(S_{k}\left(\zeta^{n}\right)\left(\tau\right)\right) & \le\bigl\Vert(\phi_{\omega}^{n})'\bigr\Vert_{X}\le\prod_{i=1}^{k}\bigl\Vert(\phi_{\omega_{i}}^{n})'\bigr\Vert_{X}\le R^{k}\prod_{i=1}^{k}\bigl\Vert(\phi_{\omega_{i}})'\bigr\Vert_{X}\\
 & \le R^{k}K^{k}\bigl\Vert(\phi_{\omega})'\bigr\Vert_{X}\le K^{k+1}R^{k}\exp\sup_{\rho\in I^{\N}\cap\left[\omega\right]}\left(S_{k}\left(\zeta\right)\left(\rho\right)\right).\end{align*}
Combining this with Assumption \ref{ass:AdditionalBoundOnPsi } we
have with $C:=K^{t\left(k+1\right)}R^{tk}M^{|\beta|}$ and $t\ge0$\begin{align*}
\exp\sup_{\tau\in{I}_{n}^{\N}\cap\left[\omega\right]}\left(S_{k}\left(t\zeta^{n}+\beta\psi^{n}\right)\left(\tau\right)\right) & \le\exp\sup_{\tau\in I_{n}^{\N}\cap\left[\omega\right]}\left(S_{k}\left(t\zeta^{n}\right)\left(\tau\right)\right)\exp\sup_{\tau\in I_{n}^{\N}\cap\left[\omega\right]}\left(S_{k}\left(\beta\psi^{n}\right)\left(\tau\right)\right)\\
 & \le C\exp\sup_{\rho\in I^{\N}\cap\left[\omega\right]}\left(S_{k}\left(t\zeta\right)\left(\rho\right)\right)\exp\inf_{\tau\in I^{\N}\cap\left[\omega\right]}\left(S_{k}\left(\beta\psi\right)\left(\tau\right)\right)\\
 & \le C\exp\sup_{\rho\in I^{\N}\cap\left[\omega\right]}\left(S_{k}\left(t\zeta+\beta\psi\right)\left(\rho\right)\right).\end{align*}
For $t<0$ we have by the assumption (\ref{enu:condition2}) with
$\widetilde{K}:=\sup K_{\Phi_{n}}<\infty$ that\begin{alignat*}{1}
\exp\sup_{\tau\in I_{n}^{\N}\cap\left[\omega\right]}\left(S_{k}\left(\zeta^{n}\right)\left(\tau\right)\right) & \ge\widetilde{K}^{-1}\bigl\Vert(\phi_{\omega}^{n})'\bigr\Vert_{X}\ge\widetilde{K}^{-k-1}\prod_{i=1}^{k}\bigl\Vert(\phi_{\omega_{i}}^{n})'\bigr\Vert_{X}\\
 & \ge R^{-k}\widetilde{K}^{-k-1}\prod_{i=1}^{k}\bigl\Vert(\phi_{\omega_{i}})'\bigr\Vert_{X}\ge R^{-k}\tilde{K}^{-k-1}\bigl\Vert(\phi_{\omega})'\bigr\Vert_{X}\\
 & \ge R^{-k}\widetilde{K}^{-k-1}\exp\sup_{\rho\in I^{\N}\cap\left[\omega\right]}\left(S_{k}\left(\zeta\right)\left(\rho\right)\right).\end{alignat*}
As above we have with $C:=R^{-tk}\widetilde{K}^{-t\left(k+1\right)}M^{|\beta|}$
and $t<0$\begin{align*}
\exp\sup_{\tau\in{I}_{n}^{\N}\cap\left[\omega\right]}\left(S_{k}\left(t\zeta^{n}+\beta\psi^{n}\right)\left(\tau\right)\right) & \le\exp\sup_{\tau\in I_{n}^{\N}\cap\left[\omega\right]}\left(S_{k}\left(t\zeta^{n}\right)\left(\tau\right)\right)\exp\sup_{\tau\in I_{n}^{\N}\cap\left[\omega\right]}\left(S_{k}\left(\beta\psi^{n}\right)\left(\tau\right)\right)\\
 & \le C\exp\sup_{\rho\in I^{\N}\cap\left[\omega\right]}\left(S_{k}\left(t\zeta\right)\left(\rho\right)\right)\exp\inf_{\tau\in I^{\N}\cap\left[\omega\right]}\left(S_{k}\left(\beta\psi\right)\left(\tau\right)\right)\\
 & \le C\exp\sup_{\rho\in I^{\N}\cap\left[\omega\right]}\left(S_{k}\left(t\zeta+\beta\psi\right)\left(\rho\right)\right).\end{align*}

\textbf{ad (\ref{enu:condition1})}: Since the potential $\psi$ is
bounded and the topological entropy infinite we have $\P\left(\beta\psi\right)=\infty$.
Hence, to verify the condition in Definition \ref{def:regularconvergence-1}
we only have to concider the case $t>0$, since for $t\le0$ we have
$\P\left(t\zeta+\beta\psi\right)\ge\P\left(\beta\psi\right)=\infty$.
But this case has been treated in \textbf{(\ref{enu:condition2})}
without any additional assumption on $\Phi$. 

\textbf{(\ref{enu:condition3})}$\pmb{\implies}$\textbf{(\ref{enu:condition2})}:
Since by Lemma \ref{lem:s_n converges to s } the contraction ratios
$s_{\Phi^{n}}$ of $\Phi^{n}$ defined in (\ref{eq:def s_phi}) converge
to $s_{\Phi}<1$ we conclude by Lemma \ref{l2.033101} that $\sup_{n}$$K_{\Phi^{n}}<\infty$. 

\textbf{(\ref{enu:condition4})}$\pmb{\implies}$\textbf{(\ref{enu:condition2})}:
This implication is an immediate consequence of \cite[Theorems 4.1.2 and 4.1.3]{urbanskimauldin-gdmsMR2003772}.
See also \emph{Proof of Claim} in the proof of Theorem 5.20 in \cite{lambdatopologypreprint}
for a similar argument.\end{proof}
\begin{acknowledgement*}
We would like to thank Mario Roy, Hiroki Sumi and Mariusz Urba{\'n}ski
for their helpful comments on an earlier draft of this paper.
\end{acknowledgement*}
\newcommand{\etalchar}[1]{$^{#1}$}


\begin{thebibliography}{Wal82}
\bibitem[FP85]{FrischParisi:85} U.~Frisch and G.~Parisi, \emph{On
the singularity structure of fully developed turbulence}, Turbulence
and predictability in geophysical fluid dynamics and climate dynamics
(North Holland Amsterdam), 1985, pp.~84--88.

\bibitem[HJK{\etalchar{+}}86]{HalseyJensenKadanoff:86} T.~C.~Halsey,
M.~H.~Jensen, L.~P.~Kadanoff, I.~Procaccia, and B.~J.~Shraiman,
\emph{Fractal measures and their singularities: The characterization
of strange sets}, Phys. Rev. A \textbf{85} (1986), no.~33, 1141--1151.

\bibitem[JK10]{JaerischKess:08} J.~Jaerisch and M.~Kesseb{\"o}hmer,
\emph{The arithmetic-geometric scaling spectrum for continued fractions},
Arkiv f{\"o}r Matematik \textbf{48} (2010), no. 2, doi:10.1007/s11512--009--0102--8.

\bibitem[JKL10]{JaerischKessLamei:10} J.~Jaerisch, M.~Kesseb{\"o}hmer,
and S.~Lamei, \emph{Induced topological pressure for countable state
{M}arkov shifts}, preprint in arXiv (2010).

\bibitem[KMS10]{KessMundayStatmann:10} M.~Kesseb{\"o}hmer, S.~Munday,
and B.~O. Stratmann, \emph{Strong renewal theorems and {L}yapunov
spectra for $\alpha$-{F}arey-{L}{\"u}roth and $\alpha$-{L{\"u}roth}
systems}, preprint in arXiv (2010).

\bibitem[KS07]{kessemultifractalsternbrocotMR2338129} M.~Kesseb{\"o}hmer
and B.~O.~Stratmann, \emph{A multifractal analysis for {S}tern-{B}rocot
intervals, continued fractions and {D}iophantine growth rates},
J. Reine Angew. Math. \textbf{605} (2007), 133--163. \MR{MR2338129}

\bibitem[KU07]{MR2343687} M.~Kesseb{\"o}hmer and M.~Urba{\'{n}}ski,
\emph{Higher-dimensional multifractal value sets for conformal infinite
graph directed {M}arkov systems}, Nonlinearity \textbf{20} (2007),
no.~8, 1969--1985. \MR{MR2343687 (2008f:37055)}

\bibitem[Kur66]{Kuratowski:66} K.~Kuratowski, \emph{Topology. {V}ol.
{I}}, New edition, revised and augmented. Translated from the French
by J. Jaworowski, Academic Press, New York, 1966. \MR{MR0217751
(36 \#840)}

\bibitem[Man74]{Mandelbrot:74} B.~B.~Mandelbrot, \emph{Intermittent
turbulence in self-similar cascades: divergence of high moments and
dimension of the carrier}, Journal of Fluid Mechanics Digital Archive
\textbf{62} (1974), no.~02, 331--358.

\bibitem[Man88]{Mandelbrot:88} \bysame, \emph{An introduction to
multifractal distribution functions}, Random fluctuations and pattern
growth ({C}argèse, 1988), NATO Adv. Sci. Inst. Ser. E Appl. Sci.,
vol. 157, Kluwer Acad. Publ., Dordrecht, 1988, pp.~279--291. \MR{MR988448}

\bibitem[MU03]{urbanskimauldin-gdmsMR2003772} D.~Mauldin and M.~Urba{\'{n}}ski,
\emph{Graph directed {M}arkov systems}, Cambridge Tracts in Mathematics,
vol. 148, Cambridge University Press, Cambridge, 2003, Geometry and
dynamics of limit sets. \MR{MR2003772 (2006e:37036)}

\bibitem[Pes97]{pesindimensiontheoryMR1489237} Ya.~B.~Pesin, \emph{Dimension
theory in dynamical systems}, Chicago Lectures in Mathematics, University
of Chicago Press, Chicago, IL, 1997, Contemporary views and applications.
\MR{MR1489237 (99b:58003)}

\bibitem[Roc70]{rockafellar-convexanalysisMR0274683} R.~T.~Rockafellar,
\emph{Convex analysis}, Princeton Mathematical Series, No. 28, Princeton
University Press, Princeton, N.J., 1970. \MR{MR0274683 (43 \#445)}

\bibitem[RSU09]{lambdatopologypreprint} M.~Roy, H.~Sumi, and M.~Urba{\'{n}}ski,
\emph{Lambda-topology versus pointwise topology}, Ergodic Theory Dynam.
Systems \textbf{29} (2009), no.~2, 685--713. \MR{MR2486790 (2010c:37050)}

\bibitem[RU05]{MR2183304} M.~Roy and M.~Urba{\'{n}}ski, \emph{Regularity
properties of {H}ausdorff dimension in infinite conformal iterated
function systems}, Ergodic Theory Dynam. Systems \textbf{25} (2005),
no.~6, 1961--1983. \MR{MR2183304 (2008c:37042)}

\bibitem[RU09]{RoyUrbanski:08} \bysame, \emph{Multifractal analysis
for conformal graph directed {M}arkov systems}, Discrete Contin.
Dyn. Syst. \textbf{25} (2009), no.~2, 627--650. \MR{MR2525196}

\bibitem[Sch99]{Schmeling:99} J.~Schmeling, \emph{On the completeness
of multifractal spectra}, Ergodic Theory Dynam. Systems \textbf{19}
(1999), no.~6, 1595--1616. \MR{2000k:37009}

\bibitem[SW77]{MR0479398} G.~Salinetti and R.~J.-B.~Wets, \emph{On
the relations between two types of convergence for convex functions},
J. Math. Anal. Appl. \textbf{60} (1977), no.~1, 211--226. \MR{MR0479398
(57 \#18828)}

\bibitem[Wal82]{walters-ergodictheoryMR648108} P.~Walters, \emph{An
introduction to ergodic theory}, Graduate Texts in Mathematics, vol.~79,
Springer-Verlag, New York, 1982. \MR{MR648108 (84e:28017)}
\end{thebibliography}
\end{document}